\documentclass[12pt]{article}
\usepackage{epsfig}
\usepackage{amsmath}
\usepackage{amsfonts}
\usepackage{amssymb}
\usepackage{color}
\usepackage{hyperref}

\title{Pluripolar hulls and fine analytic structure}

\author{Tomas Edlund and Said El Marzguioui}

\baselineskip

\newcounter{bean}

\newtheorem{theorem}{Theorem}[section]

\newtheorem{lemma}[theorem]{Lemma}

\newtheorem{definition}[theorem]{Definition}

\newcommand{\qed}{\hspace*{\fill}$\square$}

\begin{document}

\date{\today}
\maketitle \footnote{2000 Mathematics Subject Classification
31C40, 32U15, 32E20}

\begin{abstract}
We discuss the relation between pluripolar hulls and fine analytic
structure. Our main result is the following. For each non polar
subset $S$ of the complex plane $\mathbb C$ we prove that there
exists a pluripolar set $E \subset (S \times \mathbb C)$ with the
property that the pluripolar hull of $E$ relative to $\mathbb C^2$
contains no fine analytic structure and its projection onto the
first coordinate plane equals $\mathbb C$.
\end{abstract}

\begin{center}
\section{Introduction}\noindent
\end{center}Denote by $\Omega$ an open subset of $\mathbb C^n$ and let $E \subset \Omega$ be a pluripolar subset.
It might be the case that any plurisubharmonic function $u(z)$
defined in $\Omega$ that is equal to $-\infty$ on the set $E$ is
necessarily equal to $-\infty$ on a strictly larger set. For
instance, if $E$ contains a non polar proper subset of a connected
Riemann surface embedded into $\mathbb C^n$, then any
plurisubharmonic function defined in a neighborhood of the Riemann
surface which is equal to $-\infty$ on $E$ is automatically equal
to $-\infty$ on the whole Riemann surface. In order to try to
understand some aspect of the underlying mechanism of the
described "propagation" property of pluripolar sets, the
pluripolar hull of graphs $\Gamma_f(D)$ of analytic functions $f$
in a domain $D\subset \mathbb C$ has been studied in a number of
papers. (See for instance \cite{EW03}, \cite{edl}, \cite{pol} and
\cite{wie}.)

The pluripolar hull $E_{\Omega}^*$ relative to $\Omega $ of a
pluripolar set $E$ is defined as follows.

$$
E_{\Omega}^{*}=\bigcap \{z\in \Omega \ : u(z)= - \infty\},
$$
where the intersection is taken over all plurisubharmonic
functions defined in $\Omega$ which are equal to $-\infty$ on $E$.
The set $E$ is called \emph{complete pluripolar in $\Omega$} if
there exists a plurisubharmonic function on $\Omega$ which equals
$-\infty$ precisely on $E$.

As remarked above a necessary condition for a pluripolar set $E$
to satisfy $E_{\Omega}^{*} = E$ is that $E \cap A$ is polar in $A$
(or $E \cap A = A$) for all one-dimensional complex analytic
varieties $A \subset \Omega$. The fact that this is not a
sufficient condition was proved by Levenberg in \cite{le88}. By
using a refinement of Wermer's example of a polynomially convex
compact set with no analytic structure (cf. \cite{we}) Levenberg
proved that there exists a compact set $K \subset \mathbb C^2$
satisfying $K \neq K_{\mathbb C^2}^{*}$, and the intersection of
$K$ with any one dimensional analytic variety $A$ is polar in $A$.
In this example it is not clear what the pluripolar hull
$K_{\mathbb C^2}^{*}$ equals.

We will say that a set $S\subset \mathbb C^n$ contains fine
analytic structure if there exists a non constant map $\varphi:U
\rightarrow S$ from a fine domain $U \subset \mathbb C$ whose
coordinate functions are finely holomorphic in $U$ (see Definition
2.3 below). Such a map $\varphi$ will be called a \emph{fine
analytic curve}.

Motivated by recent results of J\"oricke and the first author (cf.
\cite{edl}), the following result was proved in \cite{EMW}.

\begin{theorem}\label{jan} Let $\varphi$ :\ $U$ $\longrightarrow$ $\mathbb C^n$ be a finely holomorphic map
on a fine domain $U \subset \mathbb{C}$ and let $E \subset \mathbb C^n$ be a pluripolar set. Then the following hold\\
(1) $\varphi(U)$ is a pluripolar subset of $\mathbb C^n$\\
(2) If $\varphi^{-1}(\varphi(U)\cap E )$ is a non polar subset of
$\mathbb C$ then $\varphi(U)\subset E^{*}_{\mathbb C^n}$.
\end{theorem}

In view of this result one may expect to get more information on
the pluripolar hull $E^{*}_{\mathbb C^n}$ by examining the
intersection of the pluripolar set $E$ with fine analytic curves.
Since many curves in $\mathbb C^n$ are complete pluripolar (see
\cite{edlund}) one cannot expect that $E^{*}_{\mathbb C^n}$ always
contains fine analytic structure. However if we consider the non
trivial part $E^{*}_{\mathbb C^n} \smallsetminus E$ the situation
is up to now slightly different. In fact, all examples we have
seen so far have the property that if $E_{\mathbb C^n}^{*}
\smallsetminus E$ is nonempty then for each $w \in E_{\mathbb
C^n}^{*} \smallsetminus E$ there exists a finely analytic curve
$\varphi$ contained in $E_{\mathbb C^n}^*$ which passes through
the point $w$. (i.e. $\varphi:U \rightarrow E_{\mathbb C^2}^* $ is
a finely analytic curve and $\varphi(z)=w$ for some $z \in U$). In
this paper we prove that no such conclusion holds in general. We
have the following main result.

\begin{theorem}\label{MAIN}
For each proper non polar subset $S \subset \mathbb C$ there
exists a pluripolar set $E \subset (S \times \mathbb C)$ with the
property that $E_{\mathbb C^2}^*$ contains no fine analytic
structure and the projection of $E_{\mathbb C^2}^*$ onto the first
coordinate plane equals $\mathbb C$.
\end{theorem}

The set $E$ will be a subset of a complete pluripolar set $X$
which is constructed in the same spirit as Wermer's polynomially
convex compact set without analytic structure. \\
Let us describe more precisely the content of the paper. In
Section 2 we briefly recall the construction of Wermer's set and
prove that it contains no fine analytic structure. This leads to
Theorem \ref{norm} which slightly generalizes a result in
\cite{le88}. The main result is proved in Section 3. Subsection
3.1 is devoted to construct the above mentioned set $X$ and in
Subsection 3.2 we show that $X$ contains no fine analytic
structure. In Subsection 3.3 we define the set $E$ and describe
$E^{*}_{\mathbb C^2}$. Finally, in
Section 4 we make some remarks and pose two open questions.\\
Readers who are not familiar with basic results on finely
holomorphic functions and fine potential theory are referred to \cite{Fu72} and \cite{Fu81}. \\

\noindent\textbf{Acknowledgments}. Part of this work was completed
while the first author was visiting Korteweg-de Vries Institute
for Mathematics, University of Amsterdam. He would like to thank
this institution for its hospitality. The authors thank Professor
Jan Wiegerinck for very helpful discussions.

\begin{center}
\section{Wermer's example}
\end{center}
In this Section we sketch the details of Wermer's construction
given in \cite{we}. Denote by $\mathcal{D}_r$ the open disk with
center zero and radius $r$ and by $\mathcal{C}_r$ the open
cylinder $\mathcal{D}_r \times \mathbb C$. Let $a_{1}, a_{2},....$
denote the points in the disk $\mathcal{D}_{\frac{1}{2}}$ with
rational real and imaginary part. For each $j$ we denote by
$B_{j}(z)$ the algebraic ($2$-valued) function
\begin{eqnarray*}
B_j(z)=(z-a_{1})(z-a_{2})...(z-a_{j-1}) \sqrt{(z-a_{j})}.
\end{eqnarray*}
To each $n$-tuple of positive constants $c_{1}, c_{2},..., c_{n}$
we associate the algebraic ($2^{n}$-valued) function $g_{n}(z) =
\sum_{j=1}^{n}c_{j}B_{j}(z)$. Let $\sum (c_{1},..., c_{n}), \ n=
1, 2,...$ be the subset of the Riemann surface of $g_{n}(z)$ which
lies in $\overline{\mathcal{C}_{\frac{1}{2}}}$.
\begin{lemma}\label{wermer}[\cite{we}, lemma 1] There exist positive constants $c_{1}, c_{2},...$,
with $c_{1}=\frac{1}{10}$ and $c_{n+1}\leq (\frac{1}{10})c_{n}$,
$n=1,2,...$ and a sequence of
polynomials $\{p_{n}(z,w)\}$ such that:\\
(1) $\{p_{n} = 0\} \cap \{ |z| \leq \frac{1}{2}\}= \sum
(c_{1},...,
c_{n}), \ n= 1, 2,...$\\
(2) $\{ |p_{n+1}| \leq \varepsilon _{n+1}\} \cap \{ |z| \leq
\frac{1}{2}\} \subset \{ |p_{n}| \leq \varepsilon _{n}\} \cap \{
|z| \leq \frac{1}{2}\}$, $n= 1, 2,...$\\
(3) If $ |a| \leq \frac{1}{2}$ and $ |p_{n}(a,w)|\leq \varepsilon
_{n}, \ then \ there \ is \ a \ w_{n} \ with \ p_{n}(a,w_{n})=0 \
and \ |w-w_{n}| \leq \frac{1}{n}$, $n= 1, 2,...$.
\end{lemma}

With $p_{n}$, $\varepsilon_{n}$, $n=1,2,...$ chosen as in Lemma
$\ref{wermer}$, we put
$$
Y= \bigcap_{n=1}^{\infty}[\{ |p_{n}| \leq \varepsilon _{n}\} \cap
\{ |z| \leq \frac{1}{2}\} ].
$$
Clearly, $Y$ is a compact polynomially convex subset of
$\mathbb{C}^{2}$. It was shown by Wermer that $Y$ has no analytic
structure i.e. $Y$ contains no non-constant analytic disk. In fact
he proves something stronger. The set $Y$ defined above contains
no graph of a continuous function defined on a circle in
$\mathcal{D}_{\frac{1}{2}}$ which avoids all the branch points
$\{a_i\}$. Using this observation the following lemma follows.
\begin{lemma}\label{fine}There is no fine analytic curve contained in
$Y$.
\end{lemma}

Before we prove Lemma \ref{fine} we recall the following
definition (cf. \cite{Fu81}, page 75):

\begin{definition}\label{D1} Let $U$ be a finely open set in ${\mathbb{C}}$. A function
$f$ :\ $U$ $\longrightarrow$ $\mathbb{C}$ is said to be finely
holomorphic if every point of $U$ has a compact (in the usual
topology) fine neighbourhood $K \subset U$ such that the
restriction $f\mid_{K}$ belongs to $R(K)$.
\end{definition}

Here $R(K)$ denotes the uniform closure of the algebra of all
restrictions to $K$ of rational functions on $\mathbb{C}$ with
poles off $K$.

\textit{Proof of Lemma 2.2}. Let $\varphi$ :\ $U$ $\rightarrow$
$Y$, $z\mapsto (\varphi_{1}(z), \varphi_{2}(z))$ be a fine
analytic curve contained in $Y$. If $\varphi_{1}(z)$ is constant
on $U$ then $\varphi_{2}(z)$ must also be a constant since non
constant finely holomorphic functions are finely open maps and by
the construction of the set $Y$ the fibre $Y\cap (\{z\} \times
\mathbb C)$ is a Cantor set or a finite set for any point $z \in
\overline{\mathcal{D}}_{1/2}$. Assume therefore that
$\varphi_{1}(z)$ is non-constant. In particular, there is a point
$z_0 \in U$ where the fine derivative of $\varphi_1(z)$ does not
vanish. Hence $\varphi_{1}(z)$ is one-to-one on some finely open
neighborhood $V \subset U$ of the point $z_0$. By considering the
map $z \mapsto (\varphi_1 \circ \varphi^{-1}_1(z),\varphi_2 \circ
\varphi^{-1}_1(z))$, defined on the finely open set $\varphi_1(V)$
we may assume that $\varphi$ is of the form $z \mapsto (z,g(z))$
where $g(z)=\varphi_2 \circ \varphi^{-1}_1(z)$ is finely
holomorphic in the finely open set $V'= \varphi_1(V) \subset
\mathcal{D}_{1/2}$. By Definition 2.3 there exists a compact
subset $K \subset V'$ with non-empty fine interior such that
$g(z)$ is a continuous function on $K$ (with respect to the
Euclidean topology). Shrinking $K$ if necessary we may assume that
$K\cap\{a_{1},a_{2},....\}=\emptyset$. Let $p$ be a point in the
fine interior of $K$. It is well known that there exists a
sequence of circles $\{C(p,r_j)\}$ contained in $K$ with centers
$p$ and radii $r_j \rightarrow 0$ as $j \rightarrow \infty$.
Clearly, the circle $C(p,r_j)$ avoids the branch points $\{a_{1},
a_{2},....\}$ and its image under the continuous map $z \mapsto
(z,g(z))$ is contained in $Y$. By the above observation this is
not possible. Hence $Y$
contains no fine analytic structure. \qed  \\

Denote by $d_{n}$ the degree of the one variable polynomial $w
\mapsto p_{n}(z,w)$ where $p_n(z,w)$ is the polynomial given in
Lemma 2.1. Assume that the set $Y$ is constructed using the
parameters $\epsilon_{n}$ satisfying the following condition
\begin{equation}\label{lev1}
 \lim_{n\rightarrow \infty}(\epsilon_{n})^{1/d_{n}}=0 .
\end{equation}
It is shown in \cite{le2} that with this choice the set $Y \cap
\mathcal{C}_{1/2}$ is complete pluripolar in $\mathcal{C}_{1/2}$.
Using this result and Lemma \ref{fine} we are able to generalize a
result in \cite{le88}.
\begin{theorem}\label{norm} Fix $\delta \in (0,1/2)$ and let $ Y_{\delta}= \bigcap_{n=1}^{\infty}[\{ |p_{n}|
\leq \varepsilon _{n}\} \cap \{ |z| \leq \delta \} ]$ be
constructed using the parameters $\varepsilon_{n}$ satisfying
$(\ref{lev1})$. Then \\ \noindent (a) $\varphi^{-1}(\varphi(U)\cap
Y_{\delta})$ is a polar subset of $U$ for all fine analytic
curves $\varphi$ : $U \rightarrow \mathbb{C}^{2}$. \\
\noindent(b) $Y_\delta \neq (Y_\delta)_{\mathbb C^2}^*$.
\end{theorem}

\textit{Proof of Theorem \ref{norm}}. In order to prove (a) we
argue by contradiction. Assume therefore that $\varphi:U
\rightarrow \mathbb{C}^{2}$ is a fine analytic curve and
$\varphi^{-1}(\varphi(U)\cap Y_{\delta})$ is a non polar subset of
$U$. Then there is a fine domain $U_{k_0} \subseteq U$ such that
$\varphi(U_{k_0} )\subset {\mathcal{C}_{1/2}}$ and
$\varphi^{-1}(\varphi(U_{k_0})\cap Y_{\delta})$ is non polar.
Indeed, the set $\varphi^{-1}(\varphi(U)\cap \mathcal{C}_{1/2})$
is a finely open subset of $U$ and hence has at most countably
many finely connected components $\{U_{k}\}_{k=1}^\infty$.
Moreover, $\varphi^{-1}(\varphi(U)\cap Y_{\delta} ) \cap U_{k_0}$
is non polar for some natural number $k_0$, since otherwise
$\bigcup_{k=1}^\infty \{\varphi^{-1}(\varphi(U)\cap Y_{\delta})
\cap U_{k}\} = \varphi^{-1}(\varphi(U)\cap Y_{\delta})$ would be
polar contrary to our assumption. Since $Y \cap \mathcal{C}_{1/2}$
is complete pluripolar in $\mathcal{C}_{1/2}$ there exists a
plurisubharmonic function $u$ defined in $\mathcal{C}_{1/2}$ which
is equal to $- \infty$ exactly on $Y \cap \mathcal{C}_{1/2}$. The
function $u \circ \varphi $ is either finely subharmonic on
$U_{k_0}$ or identically equal to $-\infty$ (cf. \cite{EMW}, Lemma
3.1). Since $u$ equals $- \infty$ on the non polar subset
$\varphi^{-1}(\varphi(U)\cap Y_{\delta} ) \cap U_{k_0}$, it must
be identically equal to $- \infty$ on $U_{k_0}$. Therefore
$\varphi(U_{k_0}) \subset \{u = -\infty \} =Y \cap
\mathcal{C}_{1/2}$ contradicting Lemma \ref{fine} and (a) follows.

The proof of assertion (b) follows immediately from the proof of
Proposition 3.1 in \cite{le88}. Indeed, if $u$ is a
plurisubharmonic function defined in $\mathbb C^2$ which equals $-
\infty$ on $Y_\delta$ then the function $z \mapsto
\textrm{max}\{u(z,w): (z,w) \in Y \}$ is subharmonic in
$\mathcal{D}_{1/2}$ and since it equals $-\infty$ on
$\mathcal{D}_\delta$ it equals $-\infty$ on $\mathcal{D}_{1/2}$.
Consequently $Y\cap \mathcal{C}_{1/2} \subset
(Y_{\delta})_{\mathbb C^2}^*$ and hence $Y_\delta \neq
(Y_\delta)_{\mathbb C^2}^*$. \qed \\

\textbf{Remark}. It follows from the argument used in the proof of
assertion (b) in Theorem 2.4 that $Y\cap \mathcal{C}_{1/2} \subset
(Y_\delta)_{\mathcal{C}_{1/2}}^*$. Since the first set is complete
pluripolar in $\mathcal{C}_{1/2}$ it follows that
$(Y_\delta)_{\mathcal{C}_{1/2}}^* = Y \cap \mathcal{C}_{1/2}$.
Consequently, $(Y_\delta)_{\mathcal{C}_{1/2}}^*$ contains no fine
analytic structure. It would be nice to determine what the set
$(Y_{\delta})_{\mathbb C^2}^*$ equals and to figure out whether
this set contains fine analytic structure. We are unable to do
this. But by modifying Wermer's construction, we will in the next
Section construct a complete pluripolar Wermer-like set $X\subset
\mathbb{C}^2 $ with the property that $(X\cap (S \times
\mathbb{C}))_{\mathbb C^2}^*$ contains no fine analytic structure
for all non polar subset $S \subset \mathbb{C}$.
\begin{center}
\section{Proof of Theorem 1.2}\label{s1}\noindent
\end{center}
\subsection{Construction of the set $X$}
In this Subsection we construct the set $X$. Denote by
$\{a_k\}_{k=1}^\infty$ the points in the complex plane both of
whose coordinates are rational numbers. Without loss of generality
we may assume that $a_k \in \mathcal{D}_{k}$. For any sequence of
points $\{a_l\}_{l=1}^j$ we denote by $B_j(z)$ the algebraic
function
\begin{equation*}\label{B_j}
B_j(z) = (z-a_1) \dots (z-a_{j-1}) \sqrt{(z-a_j)}.
\end{equation*}
Denote by $\gamma_j$ a simple smooth curve with endpoints $a_j$
and $\infty$. For each $j$ $B_j(z)$ has two single-valued analytic
branches on $\mathbb C \smallsetminus \gamma _j$. Following the
notation in \cite{we} we choose one of the branches $B_j(z)$
arbitrarily and denote it by $\beta_j(z)$. Then $|\beta_j(z)|
=|B_j(z)|$ is continuous on $\mathbb C$.

For each $n+1$-tuple of positive constants $(c_1,c_2, \dots
,c_{n+1})$ we denote by $g_n(z)$ the algebraic function defined
recursively in the following way. Put $g_1(z) = c_1B_1(z)$ and
$g_2(z)= c_1B_1(z) + c_2B_2(z)$ and if $g_{n}(z)$ has been chosen
we will choose $g_{n+1}(z)$ as described below. Put $Z_{1}(z) = 1$
and for $n=2,3,\dots$ define the function $Z_{n}(z)$ as follows.
Denote by $z_1,z_2, \dots, z_l$ all the zeros of all possible
different differences $h_j(z)-h_i(z)$ ($i \neq j$) of branches
$h_i(z),h_j(z)$ of the function $g_{n}(z)$. Suppose $z_k$ is a
zero of $h_j(z)-h_i(z)$ of order $m_k$ and put $Z_{n}(z)=
\Pi_{i=1}^l(z-z_i)^{m_i}$. Note that the zeros of $Z_{n}(z)$ are
also zeros of the function $Z_{n+1}(z)$ of the same or greater
multiplicity. Define $g_{n+1}(z)=g_{n}(z)
+c_{n+1}Z_{{n}}(z)B_{n+1}(z)$.

By $\Sigma(c_1, \dots, c_n)$ we mean the Riemann surface of
$g_n(z)$ which lies in $\mathbb C^2$. In other words, $\Sigma(c_1,
\dots, c_n) = \{(z,w): z \in \mathbb C, w = w_j, j=1,2, \dots, 2^n
\},$ where $w_j, \ j=1,2, \dots, 2^n$ are the values of $g_n(z)$
at $z$.

We will choose positive constants $c_n$, $\epsilon_n$ and
polynomials $p_n(z,w)$ recursively so that
\begin{eqnarray}
\label{p1}
& &\{p_n(z,w)=0 \} \cap \mathcal{C}_{n+1}  = \Sigma(c_1,c_2,\dots,c_n) \cap \mathcal{C}_{n+1} \textrm{ and} \\
\label{p2} & &\{|p_{n+1}(z,w)| \leq \epsilon_{n+1} \} \cap
\mathcal{C}_{n+1} \subset \{|p_{n}(z,w)| < \epsilon_{n} \} \cap
\mathcal{C}_{n+1}
\end{eqnarray}
hold for $n=1,2, \dots$. The set $X$ will be of the form
\begin{equation}\label{X}
X= \bigcup_{n=1}^\infty \Big( \bigcap_{j=n}^\infty \{ |p_j(z,w)|
\leq \epsilon_j \} \cap  \mathcal{C}_{n+1} \Big ).
\end{equation}

Put $c_1=1$ and let $p_1(z,w)= w^2 -(z-a_1)$. It is clear that
$\Sigma(c_1) \cap \mathcal{C}_{2} = \{p_1(z,w) = 0 \} \cap
\mathcal{C}_{2}$. Choose $\epsilon_1>0$ so that if $z_0 \in
\mathcal{D}_2$ and $|p_1(z_0,w)| \leq \epsilon_1$ then there
exists $(z_0,w_1) \in \Sigma(c_1)\cap \mathcal{C}_2$ with $|w-w_1|
\leq 1$. Let $\mathcal{B}_2 = \mathcal{D}_2 \times
\mathcal{D}_{\rho_1}$ be a bidisk where $\rho_1$ is chosen so that
\begin{displaymath}
\{ |p_1(z,w)| \leq \epsilon_1 \} \cap \mathcal{C}_2 = \{|
p_1(z,w)| \leq \epsilon_1 \} \cap  \mathcal{B}_2.
\end{displaymath}

Assume that $c_n, \epsilon_n$ and $p_n(z,w)$ have been chosen so
that (\ref{p1}) and (\ref{p2}) hold. We will now choose $c_{n+1}$
and $p_{n+1}(z,w)$.  We denote by $w_j(z)$, $j =1,2,\dots, 2^n$
the roots of $p_n(z,\cdot) = 0 $ and to each positive constant $c$
we assign a polynomial $p_c(z,w)$ by putting
\begin{eqnarray}
\label{pc1}
p_c(z,w) = \Pi_{j=1}^{2^n} \Big ( (w-w_j(z))^2 -
c^2(Z_{n}(z)B_{n+1}(z))^2 \Big ).
\end{eqnarray}
Then $p_c(z,\cdot) =0 $ has the roots $w_j(z) \pm
cZ_{n}(z)B_{n+1}(z)$, $j = 1,2, \dots ,2^n$ and so $$\{p_c(z,w) =0
\} = \Sigma (c_1,c_2,\dots,c_n,c).$$ Note that from (\ref{pc1})
$$p_c= p_n^2 + c^2q_1 + ... +(c^2)^{2^n}q_{2^n},$$ where the $q_j$
are polynomials in $z$ and $w$, not depending on $c$. Choose $c >
0$ so that

\begin{eqnarray}
\label{c1}
& &\Sigma(c_1,c_2,\dots,c_n,c) \cap \mathcal{C}_{n+1} \subset \{|p_{n}(z,w)| < \epsilon_{n}/2 \} \cap \mathcal{C}_{n+1} \textrm{ and}\\
\label{c2} & & c \cdot |Z_{n}(z)B_{n+1}(z)| \leq (1/10) c_n
|Z_{{n-1}}(z)B_n(z)| \textrm{ holds for all } z\in
\mathcal{D}_{n+1}.
\end{eqnarray}

Decreasing $c$ if necessary we may assume that if $h_i(z)$ and
$h_j(z)$ are any different branches of the function $g_n(z)$ the
estimate
\begin{equation}\label{XXX}
 |h_j(z)-h_i(z)| \geq 2c |Z_{n}(z)B_{n+1}(z)|
\end{equation}
holds in $\mathcal{D}_{n+1}$ with equality exactly at the zeros of
$Z_{n}(z)$ which are contained in $\mathcal{D}_{n+1}$ and at the
points $a_1, \dots a_n$. This estimate will be needed later when
we prove that $X$ contains no fine analytic structure. Choose
$c_{n+1} = c$.

Let $\mathcal{B}_{n+2} = \mathcal{D}_{n+2} \times
\mathcal{D}_{\rho_{n+2}}$ be a bidisk where $\rho_{n+2}$ is chosen
so that $\{|p_{n}(z,w)| \leq \epsilon_{n} \} \cap
\mathcal{C}_{n+2} = \{|p_{n}(z,w)| \leq \epsilon_{n} \} \cap
\mathcal{B}_{n+2}$ and $\rho_{n+2} > \rho_{n+1}+1$. Let $\delta
>0$ be a constant such that $|\delta \cdot p_c(z,w)| < 1$ in
$\mathcal{B}_{n+2}$ and choose $p_{n+1}(z,w) =\delta \cdot
p_c(z,w)$.

We now turn to the choice of $\epsilon_{n+1}$. Since the part of
the zero set of $p_{n+1}(z,w)$ which is contained in
$\mathcal{B}_{n+1}$ is a subset of $\{|p_n(z,w)| < \epsilon_n/2 \}
\cap \mathcal{B}_{n+1}$ it is possible to find a natural number
$m_{n+1}$ so that
\begin{equation}\label{es4}
 \frac{1}{m_{n+1}} \log |p_{n+1}(z,w)| \geq - \frac{1}{2^n} \textrm{ for all } (z,w)  \in \mathcal{B}_{n+1} \smallsetminus \{|p_{n}(z,w)| \leq \epsilon_{n} \}.
\end{equation}
Choose $\epsilon_{n+1} < \epsilon_n$ so that
\begin{equation}\label{es1}
 \frac{1}{m_{n+1}}\log |p_{n+1}(z,w)| \leq -1 \textrm{ for all }(z,w)  \in \{|p_{n+1}(z,w)| \leq \epsilon_{n+1} \} \cap \mathcal{C}_{n+2}.
\end{equation}
By decreasing $\epsilon_{n+1}$ we may assume that ($\ref{p2}$) and
the following assumption hold.
\begin{equation}\label{es11}\begin{split}
{}&\textrm{If } (z_0,w) \in \mathcal{C}_{n+2}  \textrm{ and }
|p_{n+1}(z_0,w)| \leq \epsilon_{n+1}, \textrm{ then there exists }
\\ {}& (z_0,w_n) \in \mathcal{C}_{n+2} \textrm{ such that }
|p_{n+1}(z_0,w_n)| =0 \textrm{ and } |w-w_n| \leq 1/n.
\end{split}\end{equation}
This ends the recursion.
\begin{lemma} The set $X$ defined by (\ref{X}) is complete pluripolar in
$\mathbb C^2$.
\end{lemma}

\textit{Proof.} Define for $n \geq 2$ the plurisubharmonic
function
\begin{equation*}\label{es2}
 u_n (z,w) = \max \big \{ \frac{1}{m_n}\log |p_n(z,w)|, -1 \big \}
\end{equation*}
and put $u(z,w) =\sum_{n \geq 2} u_n(z,w)$. Then $u(z,w)$ is
plurisubharmonic in $\mathbb C^2$. Indeed, since the bidisks
$\mathcal{B}_n$ exhaust $\mathbb C^2$ and $|p_{n}(z,w)| <1$ in
$\mathcal{B}_{n+1}$ the series $\sum_{n \geq 2} u_n(z,w)$ will be
decreasing on each fixed bidisk $\mathcal{B}_N$ after a finite
number of terms and hence plurisubharmonic there. Since
plurisubharmonicity is a local property $u(z,w)$ is
plurisubharmonic in $\mathbb C^2$. If $(z_0,w_0) \in X$, then for
some natural number $N$,
$ (z_0,w_0) \in \bigcap_{j=N}^{\infty} \{ |p_j(z,w)| \leq
\epsilon_j \} \cap \mathcal{C}_{N+1}. $
Condition (\ref{es1}) above implies that
$u(z_0,w_0) =Const + \sum_{n > N} u_n(z_0,w_0)  = - \infty.$
Finally if $(z_0,w_0) \notin X$ then there exists a natural number
$N$ such  that $(z_0,w_0) \in \mathcal{B}_N$ and $(z_0,w_0) \notin
\{ |p_n(z,w)| \leq \epsilon_n \} \cap \mathcal{B}_{N}$ for all $n
\geq N$. By (\ref{es4})
\begin{displaymath}
u(z,w) =Const + \sum_{n > N} \max \big \{ \frac{1}{m_n}\log
|p_n(z,w)|, -1 \big \} \geq Const +\sum_{n > N} -\frac{1}{2^n} > -
\infty.
\end{displaymath}
The Lemma follows. \qed  \\

\subsection{$X$ contains no fine analytic
structure}\label{s3}\noindent In this Section we show that $X$
contains no fine analytic structure. Suppose that $z \mapsto
(\varphi_1(z),\varphi_2(z))$ is a fine analytic curve whose image
is contained in $X$. If $\varphi_1(z)$ is constant then
$\varphi_2(z)$ must be constant since $X \cap (\{z_0\} \times
\mathbb C)$ is a Cantor set or a finite set for any point $z_0 \in
\mathbb C$. On the other hand, if $\varphi_1(z)$ is non-constant,
then using the arguments given in the proof of Lemma \ref{fine} we
may assume that the fine analytic curve contained in $X$ is given
by $z\mapsto (z,m(z))$ where $m(z)$ is a finely holomorphic
function defined in $U$ where $U \subset \mathcal{D}_n$ for some
natural number $n$. Fix a point $z' \in U \smallsetminus \{a_1,
\dots, a_n \}$ . By the definition of finely holomorphic functions
we can find a compact (in the usual topology) fine neighborhood $K
\subset U $ of $z'$ where $m(z)$ is continuous. Shrinking $K$ if
necessary we may assume that $(K \smallsetminus \{z'\})
\cap(\{a_j\}_{j=1}^\infty\cup\{Z_{{k-1}}(z)=0\}_{k=2}^\infty) =
\emptyset$. Since the complement of $K$ is thin at $z'$, one can
find a sequence of circles $\{C(z',r_i)\} \subset K$ with $r_i
\rightarrow 0$ as $i \rightarrow \infty$. Choose one of the
circles $C(z',r_j)$ so that none of the points $a_1, \dots, a_n$
are contained in $\{|z-z'| \leq r_j \}$. Let $a_k$ be the first
point in the sequence $\{a_j\}_{j=n+1}^\infty$ which is contained
in $\{|z-z'| \leq r_j \}$. Note that $a_k \in \{|z-z'| < r_j \}$,
$m(z)$ is continuous on $C(z',r_j)$ and the function $Z_{k-1}(z)
\beta_{k}(z) \neq 0$ when $z \in C(z',r_j)$. The fact that the
image of $C(z',r_j)$ under the map $z \mapsto (z,m(z))$ is a
subset of $X$ will lead us to a contradiction and hence $X$
contains no fine analytic structure. In order to prove this fix a
point $z_1 \in C(z',r_j) $ and denote by $\Re$ the $2^k$ branches
of the algebraic function $g_k(z)$ defined on $C(z',r_j)
\smallsetminus \{z_1\}$.
\begin{lemma}\label{claim4} If $h_i(z)$ and $h_j(z)$ are any different functions from
$\Re$ then
\begin{equation}\label{estest}
|h_i(z) -h_j(z)| > (3/2)c_k|Z_{{k-1}}(z) \beta_k(z)|
\end{equation}
holds for all $z \in C(z',r_j) \smallsetminus \{z_1\}$.\\
\end{lemma}

\textit{Proof.} This is follows directly from (\ref{XXX}) since
$C(z',r_j) \subset \mathcal{D}_n$ and $C(z',r_j)$ does not
intersect any of the branch points $a_1, \dots, a_k$ or the zeros
of $Z_{{k-1}}(z)$. \qed

From now on the proof that $X$ contains no fine analytic structure
follows the arguments given in \cite{we}.
\begin{lemma}\label{claim5} Fix $z_0$ in $C(z',r_j) \smallsetminus \{z_1\}$. There exists a
function $h_i(z) \in \Re$, where $h_i(z)$ depends on $z_0$ such
that
\begin{equation}\label{6}
|m(z_0)-h_i(z_0)| < (1/4)c_k|Z_{{k-1}}(z_0)\beta_k(z_0)|
\end{equation}
\end{lemma}

\textit{Proof.} By (\ref{es11}) there exists $N \geq k$ and $w_N$
such that $(z_0,w_N)$ lies on $\Sigma(c_1,\dots,c_N)$ and $m(z_0)=
w_N + R(z_0)$ where $|R(z_0)| \leq (1/10)c_k|Z_{{k-1}}(z_0)
\beta_k(z_0)|$. Thus

\begin{eqnarray*}
m(z_0)&=&   \pm c_{1} \beta_{1}(z_0) + \sum_{\nu
=2}^N \pm c_{\nu} Z_{\nu-1}(z_0)\beta_{\nu}(z_0) + R(z_0)= \\
&\overset{\text{def}}=& h_i(z_0) + \sum_{\nu =k+1}^N c_{\nu}
Z_{\nu-1}(z_0) \beta_{\nu}(z_0) +R(z_0).
\end{eqnarray*}
Since $C(z',r_j) \subset \mathcal{D}_{n+1}$ and the constants
$c_\nu$ are chosen so that (\ref{c2}) holds,
\begin{eqnarray*}
|m(z_0)-h_i(z_0)| &\leq& \sum_{\nu =k+1}^N
c_{\nu}|Z_{\nu-1}(z_0)\beta_{\nu}(z_0)| +|R(z_0) | \leq \\
&\leq& c_{k}|Z_{k-1}(z_0)\beta_{k}(z_0)| (\frac{1}{10} +
\frac{1}{10^2} + \dots ) +|R(z_0)| = \\ &=& \frac{1}{9}
c_k|Z_{{k-1}}(z_0)\beta_{k}(z_0)| + \frac{1}{10}
c_k|Z_{{k-1}}(z_0)\beta_{k}(z_0)| < \\
&<& (1/4)c_k|Z_{{k-1}}(z_0)\beta_k(z_0)|.
\end{eqnarray*}
Hence (\ref{6}) holds and the Lemma is proved. \qed
\begin{lemma}\label{claim6} Fix $z_0 \in C(z',r_j) \smallsetminus \{z_1\}$ and
let $h_i(z) \in \Re$ satisfy (\ref{6}). Then for all $z$ in
$C(z',r_j) \smallsetminus \{z_1\}$
\begin{equation}\label{t2}
|m(z)-h_i(z)| < (1/3)c_k|Z_{{k-1}}(z)\beta_{k}(z)|.
\end{equation}
\end{lemma}

\textit{Proof.} The set $\mathcal{O} = \{z \in C(z',r_j)
\smallsetminus \{z_1\}: (\ref{t2})$ holds at $z \}$ is open in
$C(z_0,r_j) \smallsetminus \{z_1\}$ and contains $z_0$. If
$\mathcal{O} \neq C(z',r_j) \smallsetminus \{z_1\}$ then there is
a boundary point $p$ of $\mathcal{O}$ on $C(z',r_j) \smallsetminus
\{z_1\}$ for which
\begin{equation}\label{t3}
|m(p)-h_i(p)| = (1/3)c_k|Z_{{k-1}}(p) \beta_{k}(p)|
\end{equation}
holds. By Lemma \ref{claim5} there is some $h_j(z)$ in $\Re$ such
that
\begin{equation}\label{t4}
|m(p)-h_j(p)| < (1/4)c_k|Z_{{k-1}}(p)\beta_{k}(p)|.
\end{equation}
Thus $|h_i(p)-h_j(p)| \leq (7/12)c_k|Z_{{k-1}}(p)\beta_k(p)|.$
Also $h_i(z) \neq h_j(z)$, in view of (\ref{t3}) and (\ref{t4}).
This contradicts Lemma \ref{claim4}. Thus $\mathcal{O} = C(z',r_j)
\smallsetminus \{z_1\}$ and Lemma \ref{claim6} follows. \qed  \\

For each continuous function $v(z)$ defined on $C(z',r_j)
\smallsetminus \{z_1\}$ which has a jump at $z_1$ we write
$L^+(v)$ and $L^-(v)$ for the two limits of $v(z)$ as $z
\rightarrow z_1$ along $C(z',r_j)$. Then, by (\ref{t2}),
\begin{displaymath}
|L^+(m)-L^+(h_i)| \leq (1/3)c_k|Z_{{k-1}}(z_1)\beta_{k}(z_1)|
\end{displaymath}
and
\begin{displaymath}
|L^-(m)-L^-(h_i)| \leq (1/3)c_k|Z_{{k-1}}(z_1)\beta_{k}(z_1)|,
\end{displaymath}
so
\begin{displaymath}
|(L^+(m)-L^+(h_i))-(L^-(m)-L^-(h_i))| \leq
(2/3)c_k|Z_{{k-1}}(z_1)\beta_k(z_1)|.
\end{displaymath}
Since $m(z)$ is continuous on $C(z',r_j)$ the jump of $h_i(z)$ at
$z_1$ is in modulus less than or equal to $(2/3)
c_k|Z_{{k-1}}(z_1)\beta_k(z_1)|\neq 0.$ But $h_i(z)$ is in $\Re$,
so its jump at $z_1$ has modulus at least
$2c_k|Z_{{k-1}}(z_1)\beta_k(z_1)|$. This is a contradiction.

\subsection{The sets $E$ and $E_{\mathbb C^2}^*$} Denote by $E$ the
pluripolar set $E = (S \times \mathbb C) \cap X$ where $S$ is a
non polar subset of $\mathbb C$. Since $X$ is complete pluripolar
in $\mathbb C^2$ it follows that $E_{\mathbb C^2}^* \subset X$. To
prove that $X \subset E_{\mathbb C^2}^*$ we argue as follows.
First we claim that the set $X$ is pseudoconcave. Indeed, by the
construction of the set $X$,
\begin{equation}\label{tgh}
\mathbb C^2 \smallsetminus X = \cup_{n=1}^\infty \{|p_n(z,w)|
> \epsilon_n \} \cap \mathcal{C}_{n+1}.
\end{equation}
By the choice of the polynomials $p_n(z,w)$ it follows that
\begin{equation*}
\{|p_n(z,w)|
> \epsilon_n \} \cap \mathcal{C}_{n+1} \subset \{|p_{n+1}(z,w)|
> \epsilon_{n+1} \} \cap \mathcal{C}_{n+2}.
\end{equation*}
Moreover, for each natural number $n$ the set $\{|p_n(z,w)|
> \epsilon_n \} \cap \mathcal{C}_{n+1}$ is a domain of holomorphy. Hence $\mathbb C^2
\smallsetminus X$ is a countable union of increasing domains of
holomorphy. By the Behnke-Stein Theorem $\mathbb C^2
\smallsetminus X$ is pseudoconvex and the claim follows.

Denote by $u(z,w)$ a globally defined plurisubharmonic function
which equals $-\infty$ on $E$. It is shown in \cite{Zlod} that the
function $z \mapsto \textrm{max}\{u(z,w): (z,w) \in X \}$ is
subharmonic in $\mathbb C$. Since the projection $S$ of $E$ onto
the first coordinate plane is non polar the function $z \mapsto
\textrm{max}\{u(z,w): (z,w) \in X \}$ will be identically equal to
$-\infty$ on $\mathbb C$ hence $u(z,w)= -\infty$ on the whole of
$X$ and consequently $E _{\mathbb C^2}^* = X$. This ends the proof
of Theorem \ref{MAIN}.
\section{Final remarks and open problems}\noindent
It follows immediately from Theorem \ref{jan} and the fact that
$X$ contains no fine analytic structure that if $\varphi:U
\rightarrow \mathbb C^2$ is a fine analytic curve, then the set
$\varphi^{-1}(\varphi(U) \cap X)$ is polar in $\mathbb C$.

Despite the result of Theorem \ref{MAIN} it should be mentioned
here that in the situation where one considers the pluripolar hull
of
the graph of a finely holomorphic function defined in a fine domain $D$, the following problem still remains open.\\

\noindent{\bf Problem 1.} Let $z\in \Gamma_f(D)_{\mathbb C^2}^*$.
Does this imply that there is a fine analytic curve contained in
$\Gamma_f(D)_{\mathbb C^2}^*$ which passes through the point
$z$?\\

It is proved in \cite{EW03} that the pluripolar hull relative to
$\mathbb C^n$ of a connected pluripolar $F_\sigma$ subset is a
connected set. It is a fairly easy exercise to show that the set
$X=E_{\mathbb C^2}^*$ in Theorem \ref{MAIN} is path connected, but
in general the pluripolar hull of a connected ($F_\sigma$)
pluripolar set is {\em{not}} path connected. Indeed, denote by
$f(z)$ an entire function of order 1/3. $f(1/z)$ has an essential
singularity at $0$ and in \cite{wie} Wiegerinck proved that the
graph $\Gamma_{f(1/z)}$ of $f(1/z)$ over $\mathbb C \smallsetminus
\{0\}$ is complete pluripolar in $\mathbb C^2$. Consequently, if
we put $E=\Gamma_{f(1/z)} \cup( \{0\} \times \mathbb C)$ then $E$
is complete pluripolar in $\mathbb C^2$ and hence $E_{\mathbb
C^2}^*=E$. Moreover $E$ is a connected $F_\sigma$ subset of
$\mathbb C^2$. By the famous Denjoy-Carleman-Ahlfors theorem (see
e.g. \cite{Ah}), entire functions of order 1/3 do not have finite
asymptotic values; i.e., there are no curves $\gamma$ ending at
infinity such that $f(z)$ approaches a finite value as $z
\rightarrow \infty$ along $\gamma$. Hence it is not possible to
find a path in $E_{\mathbb C^2}^*$ connecting a point on
$\Gamma_{f(1/z)}$ with a point in the set $\{0\} \times \mathbb
C$. In view of this remark it would be interesting to know the
answer to the following question. \\

\noindent{\bf{Problem 2.}} Is
$\Gamma_f(D)_{\mathbb C^2}^*$ path connected ? \\

Finally, we mention here again the following problem from
\cite{EMW}.\\
\textbf{Problem 3}. Let $K$ be a compact set in $\mathbb{C}^{n}$
and suppose that $\varphi^{-1}(K\cap \varphi(U))$ is a polar
subset of $U$ (or empty) for any fine analytic curve $\varphi$ :\
$U$ $\longrightarrow$ $\mathbb{C}^{n}$. Must $K$ be a pluripolar
subset of $\mathbb{C}^{n}$?

\noindent \textsc{Unit of applied statistics and mathematics,
Swedish University of Agriculture, Box 7013 SE-750 07, Uppsala,
Sweden}\\
tomas.edlund@etsm.slu.se

\noindent Kd\textsc{V Institute for Mathematics, University of
Amsterdam, Plantage Muidergracht, 24, 1018 TV, Amsterdam,
The Netherlands}\\
smarzgui@science.uva.nl


\begin{thebibliography}{999}
\bibitem{Ah}
Ahlfors, L. V., {\em Untersuchungen zur Theorie der konformen
Abbildung und der ganzen Funktionen}, Acta. Soc. Sci. Fenn., Nova
Ser. I, no. {\bf 9} (1930).

\bibitem{EW03}Edigarian, A., Wiegerinck, J., {\em{The pluripolar hull of the graph of a holomorphic function with polar singularities}},
Indiana Univ. Math. J.,  {\bf 52} (2003), 1663--1680.

\bibitem{EMW}Edigarian A., El Marzguioui, S., Wiegerinck, J., {\em{The image of a finely holomorphic map is
pluripolar}}, Preprint. arXiv: math/0701136.

\bibitem{edlund}
Edlund, T., {\em{Complete pluripolar curves and graphs}}, Ann.
Polon. Mat. {\bf 84 }(2004), 75-86.

\bibitem{edl}
Edlund, T., J\"oricke, B., {\em{The Pluripolar hull of a graph and
fine analytic continuation}}, Ark. Mat. {\bf 44 }(2006), 39-60.

\bibitem{Fu72} Fuglede, B., {\em{Finely Harmonic Functions}}, Springer Lecture Notes in Mathematics, {\bf
289}, Berlin-Heidelberg-New York, (1972).

\bibitem{Fu81} Fuglede, B., {\em{Sur les fonctions finement holomorphes}}, Ann. Inst. Fourier. {\bf 31.4} (1981), 57--88.

\bibitem{le88} Levenberg, N., {\em{On an example of Wermer}}, Ark. Mat.  {\bf 26 }(1988), 155--163.

\bibitem{le2}
Levenberg, N., Slodkowski, Z., {\em Pseudoconcave pluripolar sets
in $\mathbb C^2$}, Mat. Ann. {\bf 308 }(1998).

\bibitem{pol}
Levenberg, N., Poletsky, E., {\em Pluripolar hulls}, Mich. Math.
J. {\bf 46} (1999), 151--162.

\bibitem{Ra}
Ransford, T., {\em Potential Theory in the Complex Plane},
Cambridge University Press, (1994).

\bibitem{Zlod}
Slodkowski, Z., {\em Analytic set-valued functions and spectra},
Math. Ann. {\bf 256 }, (1981), 363-386.

\bibitem{we}
Wermer, J., {\em Polynomially convex hulls and analyticity}, Ark.
Mat. {\bf 20} (1982), 129-135.

\bibitem{wie}
Wiegerinck, J., {\em{Graphs of holomorphic functions with isolated
singularities are complete pluripolar}}, Mich. Math. J. {\bf 47}
(2000), 191--197.

\end{thebibliography}
\end{document}